\documentclass[11pt,reqno]{amsart}

\usepackage{amssymb,amsmath,amsthm,amsfonts}
\usepackage{hyperref}
\usepackage{enumitem}

\title[Cohomological Calibration and Curvature Constraints]{Cohomological Calibration and Curvature Constraints on Product Manifolds: A Topological Lower Bound}

\author{Alexander Pigazzini, Magdalena Toda}
\date{}

\theoremstyle{plain}
\newtheorem{theorem}{Theorem}[section]
\newtheorem{lemma}[theorem]{Lemma}
\newtheorem{corollary}[theorem]{Corollary}

\theoremstyle{definition}
\newtheorem{definition}[theorem]{Definition}

\newtheorem{remark}[theorem]{Remark}
\newtheorem{proposition}[theorem]{Proposition}
\newtheorem{conjecture}[theorem]{Conjecture}

\numberwithin{equation}{section}

\begin{document}

\begin{abstract}

We establish a quantitative relationship between mixed de Rham classes and the geometric complexity of metric connections with totally skew torsion on product manifolds where both factors are compact oriented surfaces. For any cohomologically calibrated connection $\nabla^C$ whose torsion $T$ has pure bidegree with respect to the product decomposition and whose harmonic projection represents a non-trivial mixed class $[\omega]$, we prove that on a non-empty open subset $\mathcal{V} \subset M$,
\[
\dim\bigl(\mathfrak{hol}_p^{\mathrm{off}}(\nabla^{C})\bigr)\;\geq\;
r^\sharp\;:=\;\operatorname{rank}_{\mathbb{R}}\bigl([\omega]_{\mathrm{mixed}}\bigr)-\dim\mathcal{K},
\]
with $\mathcal{K}$ an intrinsically defined obstruction space.  The bound is a topological invariant under metric deformations preserving the parallel‑form strata and provides an obstruction to the reduction of the holonomy along the product splitting $V_1\oplus V_2$. A counterexample shows the hypothesis is optimal. When the second factor contains a circle factor, we further show that $r^{\sharp}=1$ forces the torsion to survive dimensional reduction along it, so
that the failure of the holonomy to preserve the product splitting persists on the reduced product; the mixed rank alone cannot detect this.

\end{abstract}

\maketitle
\\
\textbf{MSC Classification 2020:} 53C05, 53C07, 53C29, 58A14, 57R19, 15A69

\section{Introduction}
\label{sec_1}

One of the broader motivations for the present work is the search for topological mechanisms that constrain geometric structures beyond the classical characteristic classes and curvature identities.  While the interaction between topology and Riemannian geometry has produced many celebrated rigidity and obstruction results, comparatively little is known
about how cohomological data can directly control the structure of metric connections with skew torsion and, in particular, their holonomy.  The viewpoint adopted here is therefore not to study a specific geometric class
of manifolds, but rather to identify a cohomological invariant whose non-triviality imposes explicit restrictions on the off-diagonal holonomy generated by a torsion connection.  From this perspective, the product geometry serves as the natural framework in which the invariant is defined, while the principal contribution is a topological obstruction linking the
mixed component of de Rham cohomology with the algebraic structure of the holonomy representation.

More precisely, let \(M=M_1\times M_2\) be a compact oriented product manifold endowed with a product metric \(g=g_1\oplus g_2\), inducing at each point an orthogonal splitting \(T_pM=V_1\oplus V_2\).  
We introduce a topological invariant, the \emph{mixed tensor rank} of a cohomology class 
\([\omega]\in H^3(M;\mathbb{R})\), and show that—after subtracting a correction term that captures the part of the class absorbed by parallel forms—it provides a lower bound for the dimension of the off‑diagonal holonomy subspace generated by any \emph{cohomologically calibrated} metric connection \(\nabla^{C}\) with totally skew torsion of pure bidegree.

A subtle but crucial point is that the subspaces spanned by a minimal harmonic decomposition 
are uniquely determined only under topological restrictions on the factors.  Consequently, our 
main result holds for \emph{admissible} product manifolds (Definition \ref{def:admissible}), i.e., where both factors are compact oriented surfaces.

For any admissible $M$, any cohomologically calibrated connection $\nabla^C$ with torsion of pure bidegree satisfies, on a non-empty open subset $\mathcal{V} \subset M$,
\[
\dim\bigl(\mathfrak{hol}_p^{\mathrm{off}}(\nabla^{C})\bigr)\;\geq\;
r^\sharp\;:=\;\operatorname{rank}_{\mathbb{R}}\bigl([\omega]_{\mathrm{mixed}}\bigr)-\dim\mathcal{K},
\]
where the obstruction space $\mathcal{K}$ is intrinsically defined via intersections of the mixed 
factor spaces $\mathcal{V}_{2,1}$, $\mathcal{V}_{1,2}$ with the parallel‑form strata 
$\mathcal{P}_k(M_i)$.  The bound is invariant under metric deformations preserving these strata.

When $r^\sharp>0$ the product splitting $V_1\oplus V_2$ is not invariant under the holonomy; equivalently, $\nabla^C$ does not reduce to a product connection along $V_1\oplus V_2$. The examples in Section \ref{sec:examples} illustrate both a non-trivial bound ($\Sigma_g\times T^2$ with a pure bidegree class) and a trivial bound ($T^2\times T^2$, where the non-pure bidegree and full parallelism of the torsion yield $r^\sharp=0$). Remark \ref{rem:counterexample} shows by explicit counterexample that without the admissibility hypothesis the bound fails to be well-defined.

The pure bidegree condition on the torsion---namely, that $T$ lies entirely in $\Gamma(\Lambda^2 V_1^*\otimes V_2^*)$ or entirely in 
$\Gamma(V_1^*\otimes\Lambda^2 V_2^*)$---admits a transparent 
geometric interpretation in the product setting: since 
$\Lambda^2 V_i^*$ is one-dimensional for a surface, a $3$-form of pure bidegree $(2,1)$ necessarily has the form 
$\operatorname{vol}_{M_1}\wedge\beta$ for a $1$-form $\beta$ on $M_2$, so the torsion couples the two factors through the area form of exactly one of them.  This condition arises naturally in Kaluza--Klein constructions on products of surfaces, where the torsion $3$-form is built as the wedge product of the curvature $2$-form of one factor with a connection $1$-form on the other: since the curvature of a surface is necessarily proportional to its area form, the resulting torsion has pure bidegree by construction.  We nevertheless conjecture (Conjecture \ref{conj:mixed}) that the theorem extends to torsion forms of mixed bidegree, i.e.\ when both $T^{2,1}$ and $T^{1,2}$ are simultaneously non-vanishing.

A further consequence concerns dimensional reduction.  When the second factor contains a circle factor, $M_2=N\times S^1$, Proposition \ref{prop:reduction} shows that $r^{\sharp}=1$ forces the torsion $3$-form to be non-zero on some
slice $M_1\times N$, for \emph{every} torsion representing the class, so that the reduced connection still fails to preserve the product splitting.  Since the mixed rank equals $1$ for every non-trivial class in this setting, the discrimination is carried entirely by the correction term $\dim\mathcal{K}$:
this is the precise sense in which $r^{\sharp}$, rather than
$\operatorname{rank}_{\mathbb{R}}([\omega]_{\mathrm{mixed}})$, is the operative invariant.

The point of view described above also places the present work within the broader program of understanding how global analytical and cohomological techniques can detect phenomena that are not immediately visible from local curvature computations alone.  The proof combines Hodge theory, the
K\"unneth decomposition of cohomology, and holonomy methods to derive a lower bound that depends only on the mixed K\"unneth class and on the parallel-form strata of the factor metrics, and is therefore independent of any particular coordinate representation or local geometric construction.
Although our results are established for admissible products
(Definition \ref{def:admissible}), i.e.\ products of compact oriented surfaces, the methods suggest a more general framework in which cohomological invariants may serve as computable obstructions for classes of metric connections with torsion on broader families of manifolds; a first step in this direction is Corollary \ref{cor:extended}.  In this
sense, the results presented here should be viewed less as isolated estimates for a particular geometric setting and more as evidence that topological information can impose explicit algebraic constraints on holonomy in a way that, to the best of our knowledge, has not been systematically exploited in this setting.

\section{Cohomologically Calibrated Metric Connections with Skew Torsion}

Throughout, \( \nabla^{LC} \) denotes the Levi-Civita connection of \(g\). We consider metric connections \( \nabla^C \) with totally skew torsion \(T\), i.e.
\begin{equation}
\nabla^C g=0,
\end{equation}
\begin{equation}
T(X,Y,Z):=g(T(X,Y),Z)\in\Omega^3(M),
\end{equation}
and write
\begin{equation}
\nabla^C=\nabla^{LC}+K,   
\end{equation}
\begin{equation}
K_{XYZ}=\tfrac12\,T_{XYZ}.
\end{equation}

Thus \(T^\flat:=T\) is a 3-form. The curvature \(R^C\) of \( \nabla^C \) can be expressed in terms of \(R^{LC}\), \(\nabla^{LC}T\), and quadratic torsion terms; in abstract index notation,
\begin{equation}\label{eq:curv-torsion}
R^C_{XY}Z
=
R^{LC}_{XY}Z
+ \tfrac12 \big( (\nabla^{LC}_X T)(Y,Z,\cdot)^\sharp - (\nabla^{LC}_Y T)(X,Z,\cdot)^\sharp \big)
+ \tfrac14\,Q_T(X,Y)Z,
\end{equation}
where \(Q_T\) is bilinear and quadratic in \(T\) (see, e.g., \cite{B (1987)} or standard torsion-curvature formulas). We shall use \eqref{eq:curv-torsion} only qualitatively: mixed components of \(T\) produce mixed components of \(R^C\) through both the \(\nabla^{LC}T\) and \(T\ast T\) terms.

The product structure $TM=V_1\oplus V_2$ induces a bigrading
on $\Omega^3(M)$: every $3$-form decomposes pointwise as
$T=T^{2,1}+T^{1,2}$, where
$T^{2,1}\in\Gamma(\Lambda^2V_1^*\otimes V_2^*)$ and
$T^{1,2}\in\Gamma(V_1^*\otimes\Lambda^2V_2^*)$
(components of type $(3,0)$ and $(0,3)$ vanish since $\dim V_i=2$). We say that $T$ has \emph{pure bidegree} if one of $T^{2,1}$, $T^{1,2}$ vanishes identically.
Since $g$ is a product metric, $\Delta_M=\Delta_{M_1}+\Delta_{M_2}$ preserves
the K\"unneth bigrading of $\Omega^{\bullet}(M)$; hence the Hodge decomposition is compatible with the bigrading, and the harmonic projection of a form of pure bidegree has the same pure bidegree.

\begin{definition}[Metric cohomological calibration]\label{def:calibration}
A metric connection \( \nabla^C \) with totally skew torsion \(T\), where the torsion 3-form represents a non-trivial mixed cohomology class \([\omega]\) via its harmonic projection (\(T^\flat=\omega_h +d\alpha+\delta \beta\), i.e. \([\omega]:=[\omega_h]\)), is called \emph{cohomologically calibrated}\footnote{Here the term ``calibrated'' is used in a purely cohomological sense: the torsion is constrained so that its harmonic projection represents a non-trivial mixed K\"unneth class. It is unrelated to the calibrated geometries of Harvey--Lawson (a closed form bounding volumes); no calibration inequality is involved.}. The class is \emph{mixed} if it lies outside the natural image of \(H^3(M_1)\oplus H^3(M_2)\) under the K\"unneth isomorphism:

\begin{equation}
H^3(M)\cong \bigoplus_{p+q=3} H^p(M_1)\otimes H^q(M_2).
\end{equation}
\end{definition}

Since \(M\) is compact and oriented, Hodge theory applies to \(g\): each cohomology class has a unique harmonic representative. We denote by \(\omega_h\) the harmonic representative of \([\omega]\) with respect to \(g\).
\\
For a thorough treatment of Hodge theory and cohomology on complex manifolds, see \cite{G-H (1978)}.

\subsection{Intrinsic Definition of Mixed Cohomology Subspaces}
\label{subsec:intrinsic}

We now address a subtle but fundamental point: the subspaces spanned by the factors in a minimal harmonic decomposition are not uniquely determined in general. To obtain a well--defined geometric invariant we must impose conditions on the topology of the factors.

\newcommand{\minrk}{\operatorname{mrk}}

\begin{definition}[Admissible Product Manifolds]\label{def:admissible}
A compact oriented product manifold $M = M_1 \times M_2$ is called admissible if it satisfies the following:
\begin{itemize}
\item[] both factors $M_1$, $M_2$ are compact connected oriented surfaces.
\end{itemize}
\end{definition}

\begin{lemma}[Uniqueness of Mixed Subspaces]\label{lem:invariance}
Let $M$ be an admissible product manifold and $[\omega]\in H^3(M;\mathbb{R})$
a non-trivial mixed class.  Then the subspaces
\begin{align}\label{eq:V-spaces}
\mathcal{V}_{2,1}&\subset\mathcal{H}^2(M_1)\otimes\mathcal{H}^1(M_2),\\
\mathcal{V}_{1,2}&\subset\mathcal{H}^1(M_1)\otimes\mathcal{H}^2(M_2)
\end{align}
spanned by any minimal harmonic decomposition of the mixed K\"unneth components of $[\omega]$ are \emph{intrinsically defined}, i.e. they depend only on $[\omega]$ and not on the particular minimal decomposition chosen.
\end{lemma}

\begin{proof}

For a compact oriented surface $F^2$ we have
$b_2(F)=1$ and $b_1(F)=2g(F)$.  Let $\operatorname{vol}_{F}$ be the unique
(up to scale) $g_F$--harmonic $2$--form.  Then
\[
\mathcal{H}^2(F)=\mathbb{R}\cdot\operatorname{vol}_{F}.
\]
Consequently any non-zero $T\in\mathcal{H}^2(M_1)\otimes\mathcal{H}^1(M_2)$
has rank one and can be written uniquely as
\[
T=\operatorname{vol}_{M_1}\otimes\beta,
\qquad \beta\in\mathcal{H}^1(M_2).
\]
Indeed, if $\operatorname{vol}_{M_1}\otimes\beta
=\operatorname{vol}_{M_1}\otimes\beta'$, then
$\operatorname{vol}_{M_1}\otimes(\beta-\beta')=0$ and, because
$\operatorname{vol}_{M_1}$ is non-zero as an element of $\mathcal{H}^2(M_1)$,
we must have $\beta=\beta'$.  Thus the one--dimensional subspace
$\operatorname{span}\{\beta\}$ is uniquely determined by $T$, and
\[
\mathcal{V}_{2,1}= \mathcal{H}^2(M_1)\otimes\operatorname{span}\{\beta\}
\]
is intrinsically defined.  The same argument applied symmetrically yields
the intrinsic definition of $\mathcal{V}_{1,2}$.
\end{proof}

\begin{remark}[Failure of invariance in general]\label{rem:counterexample}
Invariance fails without the admissibility hypothesis.  Let
$M_1=M_2=T^{2}\times S^{1}$, so that $M_i\cong T^{3}$ and
$b_2(M_1)=b_1(M_2)=3$, and endow each factor with the product metric $g_i=h_i\oplus d\theta_i^{2}$, where $h_i$ is a \emph{non-flat} metric on $T^{2}$.  Since a non-zero parallel $1$-form on a surface forces the Gaussian curvature to vanish identically, $\mathcal{P}_1(T^{2},h_i)=\{0\}$; by the
de Rham decomposition of parallel forms on a Riemannian product,
\[
\mathcal{P}_2(M_1)=\mathbb{R}\cdot a,\qquad
\mathcal{P}_1(M_2)=\mathbb{R}\cdot p,\qquad
\mathcal{P}_2(M_1)\otimes\mathcal{P}_1(M_2)=\mathbb{R}\,(a\otimes p),
\]
where $a:=\operatorname{vol}_{h_1}$ and $p:=d\theta_2$.  Choose
$\gamma\in\mathcal{H}^{1}(T^{2},h_1)$ and $\gamma'\in\mathcal{H}^{1}(T^{2},h_2)$
non-zero and set $b:=\gamma\wedge d\theta_1\in\mathcal{H}^{2}(M_1)$, $q:=\gamma'\in\mathcal{H}^{1}(M_2)$, so that $b\notin\mathcal{P}_2(M_1)$ and $q\notin\mathcal{P}_1(M_2)$.  The mixed class
\[
[\omega]=a\otimes p+b\otimes q\ \in\ H^{2}(M_1)\otimes H^{1}(M_2)
\]
has rank $2$ and admits the two minimal decompositions
\[
a\otimes p+b\otimes q
=\tfrac12\big[(a+b)\otimes(p+q)+(a-b)\otimes(p-q)\big],
\]
with spans $\mathcal{V}_{2,1}=\operatorname{span}\{a\otimes p,\,b\otimes q\}$ and $\mathcal{V}_{2,1}'=\operatorname{span}\{(a+b)\otimes(p+q),\,(a-b)\otimes(p-q)\}$.
Now $a\otimes p\in\mathcal{V}_{2,1}$, whereas
\[
s\,(a+b)\otimes(p+q)+t\,(a-b)\otimes(p-q)
=(s+t)\,a\otimes p+(s-t)\,a\otimes q+(s-t)\,b\otimes p+(s+t)\,b\otimes q
\]
is proportional to $a\otimes p$ only for $s=t=0$.  Hence
\[
\dim\big(\mathcal{V}_{2,1}\cap(\mathcal{P}_2(M_1)\otimes\mathcal{P}_1(M_2))\big)=1,
\qquad
\dim\big(\mathcal{V}_{2,1}'\cap(\mathcal{P}_2(M_1)\otimes\mathcal{P}_1(M_2))\big)=0,
\]
so $\dim\mathcal{K}$, and with it $r^{\sharp}$, depends on the chosen minimal decomposition: the invariant is not well defined outside the admissible class.
\end{remark}

\subsection{Mixed obstruction spaces}

Having established the intrinsic nature of the mixed factor spaces in Lemma \ref{lem:invariance} 
for admissible product manifolds (Definition \ref{def:admissible}), we now define the obstruction kernels.

\begin{definition}[Mixed obstruction spaces]\label{def:obstruction}
For an admissible product manifold $M$ and a mixed class $[\omega]$ we define
\begin{align*}
\mathcal{V}_{2,1}&:=\operatorname{span}\{\alpha_i\otimes\beta_i\}\subset
\mathcal{H}^2(M_1)\otimes\mathcal{H}^1(M_2),\\[2pt]
\mathcal{V}_{1,2}&:=\operatorname{span}\{\tilde\alpha_j\otimes\tilde\beta_j\}\subset
\mathcal{H}^1(M_1)\otimes\mathcal{H}^2(M_2),
\end{align*}
where the spans are taken over any minimal harmonic decomposition of the mixed Künneth 
components of $[\omega]$; by Lemma \ref{lem:invariance} these spans are independent of the 
particular choice.

The \emph{mixed obstruction spaces} are the intersections
\begin{align*}
\ker\Psi_p&:=\mathcal{V}_{2,1}\cap\bigl(\mathcal{P}_2(M_1)\otimes\mathcal{P}_1(M_2)\bigr),\\[2pt]
\ker\widetilde{\Psi}_p&:=\mathcal{V}_{1,2}\cap\bigl(\mathcal{P}_1(M_1)\otimes\mathcal{P}_2(M_2)\bigr),
\end{align*}
and we set \(\mathcal{K}:=\ker\Psi_p\oplus\ker\widetilde{\Psi}_p\).
The notation $\ker\Psi_p$ is anticipatory (see Remark \ref{rem:anticipatory_notation}); the intersection itself is independent of the point $p\in M$.
\end{definition}

\begin{remark}[Anticipatory notation]
\label{rem:anticipatory_notation}
The notation ``$\ker \Psi_p$'' is anticipatory: in the proof of Theorem \ref{thm:main}, we construct linear maps $\Psi_p : V_{2,1} \to T_p^*M\otimes\Lambda^3 T_p^*M$ and
$\widetilde{\Psi}_p : V_{1,2} \to T_p^*M\otimes\Lambda^3 T_p^*M$ sending each mixed simple tensor $\alpha\otimes\beta$ to the covariant derivative $\nabla^{LC}(\alpha\wedge\beta)$ of the corresponding harmonic $3$-form. A key step is to verify that their kernels coincide, on a non-empty open subset, with the intersections defined above; the induced quotient maps are then injective there, which computes $\dim\mathcal{K}$. The lower bound itself is established separately, in \emph{Case 1} of the proof of Theorem \ref{thm:main}, by exhibiting a purely off-diagonal curvature operator that lies in the holonomy algebra.
\end{remark}

\section{Off-Diagonal Curvature and Holonomy Complexity}

Let \(P_i:T_pM\to V_i\) be the orthogonal projections of the product splitting.

\begin{definition}[Off-diagonal curvature]\label{def:off-R}
For a curvature tensor \(R\), define its off-diagonal component by
\begin{equation}
R_{\mathrm{off}}(X,Y)Z := P_1 R(X,Y)(P_2 Z) + P_2 R(X,Y)(P_1 Z).
\end{equation}
\end{definition}

\begin{definition}[Off-diagonal holonomy component] \label{def:off-hol}
Let $(M,g)$ be a Riemannian product manifold with orthogonal splitting $T_pM = V_1 \oplus V_2$, and let $\nabla^C$ be a metric connection on $TM$ with totally skew-symmetric torsion. 

Let $P_i : T_pM \to V_i$ be the orthogonal projections and let $\Pi_{\mathrm{off}} : \mathfrak{so}(T_pM) \to \mathfrak{so}(T_pM)$ be the linear projection onto off-diagonal endomorphisms:
\[
\Pi_{\mathrm{off}}(A)(Z) := P_1 A(P_2 Z) + P_2 A(P_1 Z).
\]

Define the space of off-diagonal skew-symmetric endomorphisms as the image of this projection:
\[
\mathfrak{so}^{\mathrm{off}}(T_pM) := \operatorname{Im}(\Pi_{\mathrm{off}}) = \{ A \in \mathfrak{so}(T_pM) \mid A(V_1) \subset V_2, \, A(V_2) \subset V_1 \}.
\]

The \emph{off-diagonal holonomy component} is the intersection of the holonomy algebra with this space:
\[
\mathfrak{hol}_p^{\mathrm{off}}(\nabla^C) := \mathfrak{hol}_p(\nabla^C) \cap \mathfrak{so}^{\mathrm{off}}(T_pM).
\]

Equivalently, $\mathfrak{hol}_p^{\mathrm{off}}(\nabla^C) = \{ A \in \mathfrak{hol}_p(\nabla^C) \mid \Pi_{\mathrm{off}}(A) = A \}$ consists of those holonomy elements that are purely off-diagonal.
\end{definition}

\begin{remark}
Although $\mathfrak{hol}_p^{\mathrm{off}}(\nabla^C)$ is a vector subspace of $\mathfrak{hol}_p(\nabla^C)$, it is not closed under the Lie bracket, since $[\mathfrak{so}^{\mathrm{off}}, \mathfrak{so}^{\mathrm{off}}] \subset \mathfrak{so}(V_1) \oplus \mathfrak{so}(V_2)$. Nevertheless, $\dim \mathfrak{hol}_p^{\mathrm{off}}(\nabla^C) > 0$ implies that the product splitting $V_1 \oplus V_2$ is not invariant under the holonomy group: the holonomy does not act block-diagonally with respect to $V_1\oplus V_2$, so $\nabla^C$ does not reduce to a product connection along this splitting. This does \emph{not} assert that the holonomy representation is irreducible in the usual sense; other invariant subspaces, not aligned with $V_1\oplus V_2$, may still exist.
\end{remark}

By the Ambrose--Singer theorem \cite{A-S (1953)}, the holonomy algebra $\mathfrak{hol}_p(\nabla^C)$ is generated by the curvature operators of $\nabla^C$ at the points joined to $p$ by horizontal paths, hence it is independent of the choice of local generators and of the pointwise rank of the curvature tensor. The subspace $\mathfrak{so}^{\mathrm{off}}(T_pM)$, however, depends on $p$, since the splitting $V_1\oplus V_2$ is not $\nabla^{C}$-parallel; accordingly $\dim\mathfrak{hol}^{\mathrm{off}}_p(\nabla^C)$
may vary with $p$, which is why Theorem \ref{thm:main} is stated on an open subset. In particular, at each point $p\in M$ one has, clearly, $\dim\mathfrak{hol}_p(\nabla^C)\ge \dim\mathfrak{hol}_p^{\mathrm{off}}(\nabla^C)$. Because $\nabla^C$ is metric, its holonomy algebra sits in $\mathfrak{so}(T_pM)$. If $V_1\oplus V_2$ were invariant under the holonomy, all holonomy endomorphisms would be block-diagonal with respect to it, forcing $\mathfrak{hol}_p^{\mathrm{off}}(\nabla^C)=\{0\}$. Hence non-trivial off-diagonal holonomy implies that $V_1\oplus V_2$ is not holonomy-invariant.

\section{Mixed Rank of a Degree-\texorpdfstring{$3$}{3} Class}

Write the K\"unneth decomposition of \(H^3(M)\) as:
\begin{equation}
\begin{aligned}
\big(H^3(M_1)\otimes H^0(M_2)\big)\ \oplus\ \big(H^2(M_1)\otimes H^1(M_2)\big)\ \oplus\ \\
\big(H^1(M_1)\otimes H^2(M_2)\big)\ \oplus\ \big(H^0(M_1)\otimes H^3(M_2)\big).
\end{aligned}
\end{equation}

Mixed components occur precisely in bidegrees \((2,1)\) and \((1,2)\).

\begin{definition}[Mixed tensor rank]\label{def:mixed-rank}
Let \([\omega]\in H^3(M;\mathbb{R})\). Denote by \(\pi_{2,1}([\omega])\in H^2(M_1)\otimes H^1(M_2)\) and \(\pi_{1,2}([\omega])\in H^1(M_1)\otimes H^2(M_2)\) its mixed K\"unneth projections. Define
\begin{equation}
r_{2,1} = \min\Big\{ r \,:\, \pi_{2,1}([\omega])=\sum_{i=1}^r \alpha_i\otimes \beta_i,\ \alpha_i\in H^2(M_1),\ \beta_i\in H^1(M_2)\Big\},
\end{equation}
and analogously \(r_{1,2}\) for \(\pi_{1,2}([\omega])\). The mixed rank is
\begin{equation}
\operatorname{rank}_\mathbb{R}\big([\omega]_{\mathrm{mixed}}\big):= r_{2,1}+r_{1,2}.
\end{equation}
\end{definition}

For degree \(3\), the tensor rank in each bidegree is well-defined and finite, and can be computed from any decomposition of the harmonic representative \(\omega_h\) into wedge products of harmonic forms on the factors for the product metric \(g\). Although \(\omega_h\) depends on \(g\), the ranks \(r_{2,1}, r_{1,2}\) are purely topological (indeed, the tensor rank of an element in $H^{2}(M_{1})\otimes H^{1}(M_{2})$ is defined as the minimal number of simple tensors in any decomposition, hence is invariant under change of metric within the product class).

\section{Main Result and Proof}

\begin{lemma}[Existence of points of maximal rank for covariant derivatives]\label{lem:maximal-rank-open}
Let \((M_1,g_1)\) be a compact oriented Riemannian manifold and let \(\mathcal{H}^2(M_1,g_1)\) denote the space of harmonic \(2\)-forms.  
Given a subspace \(U\subset\mathcal{H}^2(M_1,g_1)\) of dimension \(r\) and its parallel part \(U_\parallel:=U\cap\mathcal{P}_2(M_1,g_1)\) of dimension \(s\), choose a basis \(\alpha_1,\dots,\alpha_r\) of \(U\) such that \(\alpha_1,\dots,\alpha_s\) span \(U_\parallel\).

If $r - s \leq 1$ (condition satisfied by admissible product manifolds: $U$ is the span of a minimal harmonic decomposition, 
so $\dim U = r_{2,1} \leq 1$, hence $r - s \leq 1$; see Definition \ref{def:mixed-rank}, Definition \ref{def:admissible} and Lemma \ref{lem:invariance}), then there exists a non-empty open set \(\mathcal{U}\subset M_1\) on which the covariant derivatives 
\[
\bigl(\nabla^{g_1}\alpha_{s+1},\dots,\nabla^{g_1}\alpha_r\bigr)
\]
are pointwise linearly independent; i.e.
\[
\mathcal{U}=\Bigl\{p\in M_1\;\Big|\;
\operatorname{rank}\bigl[(\nabla^{g_1}\alpha_i)(p)\bigr]_{i=s+1,\dots,r}=r-s\Bigr\}\neq\emptyset.
\]
Moreover, if $g_1$ is real-analytic, then $\mathcal{U}$ is dense in $M_1$.
\end{lemma}

\begin{proof}
We proceed in three steps.

\noindent\textit{1.  Reduction to a real‑analytic metric.}
Real‑analytic Riemannian metrics are dense in the space of smooth metrics on a compact manifold; this follows from the theorem of Greene--Jacobowitz on analytic approximation of smooth maps (see \cite{GreeneJacobowitz}).  Hence we can choose a sequence of real‑analytic metrics \(\{g_1^{(k)}\}_{k\in\mathbb N}\) converging to \(g_1\) in the \(C^\infty\) topology.  
The dimension of the space of parallel \(2\)-forms is upper‑semicontinuous with respect to the metric (see \cite{B (1987)}); therefore, after passing to a subsequence, we may assume
\[
\dim\mathcal{P}_2\bigl(M_1,g_1^{(k)}\bigr)\le\dim\mathcal{P}_2(M_1,g_1)=s\qquad\text{for all }k .
\]

For each \(k\) let \(\alpha_i^{(k)}\in\mathcal{H}^2\bigl(M_1,g_1^{(k)}\bigr)\) be the unique harmonic representative of the cohomology class \([\alpha_i]\).  By the stability theorem for harmonic forms under metric deformation (Warner \cite{W (1983)}, see also Kodaira--Spencer \cite{K-S (1960)}), \(\alpha_i^{(k)}\to\alpha_i\) in \(C^\infty\).  Consequently the covariant derivatives converge in \(C^1\):
\[
\nabla^{g_1^{(k)}}\alpha_i^{(k)}\longrightarrow\nabla^{g_1}\alpha_i\qquad (i=s+1,\dots,r).
\]

If we prove that for each \(k\) there exists a dense open set \(\mathcal{U}^{(k)}\subset M_1\) where the \(\nabla^{g_1^{(k)}}\alpha_i^{(k)}\) are linearly independent, then the \(C^1\) convergence guarantees that, for large enough \(k\), the same linear independence holds on a non-empty open set for the original metric \(g_1\).  Indeed, linear independence is equivalent to the non‑vanishing of at least one \((r-s)\times(r-s)\) minor of the matrix of components; each minor depends continuously on the metric and the forms, so a minor that is non‑zero on an open set for \(g_1^{(k)}\) remains non‑zero on an open set for \(g_1\) when \(k\) is sufficiently large.

We now prove the stronger density statement under the assumption that $g_1$ is real-analytic.

\noindent \textit{2.   Admissible manifolds.}
Assume now that $g_1$ is real‑analytic.  Then every harmonic form is real‑analytic (by elliptic regularity for the Hodge–Laplacian, see \cite{MorreyNirenberg}), and consequently the covariant derivatives $\nabla^{g_1}\alpha_i$ are real‑analytic sections of the vector bundle $E:=T^*M_1\otimes\Lambda^2T^*M_1$.

Consider the real‑analytic section
\[
\Phi:M_1\longrightarrow\Lambda^{r-s}E,\qquad
\Phi(p):=(\nabla^{g_1}\alpha_{s+1})(p)\wedge\cdots\wedge(\nabla^{g_1}\alpha_r)(p).
\]
The set $\mathcal{U}:=\{p\in M_1\mid\Phi(p)\neq0\}$ is exactly the set where the covariant derivatives are linearly independent; because $\Phi$ is continuous, $\mathcal{U}$ is open.  We must show that $\mathcal{U}\neq\emptyset$.

Suppose, for contradiction, that $\mathcal{U}=\emptyset$; i.e. $\Phi(p)=0$ for every $p\in M_1$.  Then at each point the sections $\nabla^{g_1}\alpha_{s+1},\dots,\nabla^{g_1}\alpha_r$ are linearly dependent.  Since the condition of rank $\le d$ is given by the vanishing of all $(d+1)\times(d+1)$ minors, which are analytic functions, the locus where the rank is constant is an analytic subset.  Let $d<r-s$ be the minimal rank and let $p_0$ be a point where the rank equals $d$.  There exists a connected open neighbourhood $V$ of $p_0$ on which the rank is constant and equal to $d$.

On $V$ we can therefore find real‑analytic functions $c_1,\dots,c_{r-s}:V\to\mathbb R$, not all zero, such that
\begin{equation}\label{eq:local-dependence}
\sum_{i=1}^{r-s}c_i(p)\bigl(\nabla^{g_1}\alpha_{s+i}\bigr)(p)=0\qquad\forall p\in V.
\end{equation}

Since we have to consider admissible manifolds (Definition \ref{def:admissible} and Lemma \ref{lem:invariance}) we have two cases:

\textit{Case $r=s$}: There are no non‑parallel forms to consider; the lemma holds vacuously.

\textit{Case $r=s+1$}: There exists a unique non‑parallel form $\alpha$ and a single analytic function $c$ such that \eqref{eq:local-dependence} reads $c(p)(\nabla^{g_1}\alpha)_p = 0$ for all $p\in V$.  If $c(p_0)\neq 0$ for some $p_0\in V$, then by analyticity $c(p)\neq 0$ on an open subset $V'\subset V$, whence $(\nabla^{g_1}\alpha)_p=0$ on $V'$. By the unique‑continuation principle for analytic sections (Aronszajn \cite{Aronszajn1957}; for differential forms see \cite{AKS1962}), $\nabla^{g_1}\alpha\equiv 0$ on the whole connected component containing $V'$, and hence on all of $M_1$. This contradicts the assumption that $\alpha$ is non‑parallel modulo $\mathcal{P}_2(M_1,g_1)$. Therefore $c\equiv 0$ on $V$, contradicting the choice of the $c_i$'s.  Hence $\mathcal{U}\neq\emptyset$. Moreover, the complement $M_1\setminus\mathcal{U} = \{p \in M_1 : (\nabla^{g_1}\alpha)(p) = 0\}$ has empty interior: if it contained an open set $V'$, then by unique continuation $\nabla^{g_1}\alpha \equiv 0$ on all of $M_1$, contradicting the hypothesis that $\alpha$ is non-parallel. Therefore $\mathcal{U}$ is dense in $M_1$.

\noindent\textit{3.  Returning to the original smooth metric.}
For the original smooth metric $g_1$, the non-emptiness
of $\mathcal{U}$ is immediate: since $\alpha_{s+1} \notin
\mathcal P_k(M_1, g_1)$, the section $\nabla^{g_1}\alpha_{s+1}$ is not
identically zero, and $\mathcal{U} = \{p \in M_1 :
(\nabla^{g_1}\alpha_{s+1})(p) \neq 0\}$ is non-empty and open.
Steps 1--2 establish the stronger conclusion that $\mathcal{U}$ is
dense when $g_1$ is real-analytic. 
\end{proof}

The non-emptiness argument holds for harmonic forms of any
degree $k$ and any smooth metric. The density statement for
real-analytic metrics likewise holds for any degree $k$, as
Steps 1--2 use only analytic approximation, semicontinuity of the dimension of parallel forms, $C^1$ convergence, and the identity principle for real-analytic sections.  Therefore, the same results hold for $H^1(M_1, g_1)$, as well as for $H^1(M_2, g_2)$ and $H^2(M_2, g_2)$.

\medskip

\begin{theorem}[Topological Lower Bound for Off-Diagonal Curvature]
\label{thm:main}

Let $M = M_1 \times M_2$ be an admissible compact oriented Riemannian product manifold 
(Definition \ref{def:admissible}) with product metric $g = g_1 \oplus g_2$.  
Let $\nabla^C$ be a cohomologically calibrated metric connection (Definition \ref{def:calibration}) whose torsion $T$ has pure bidegree with respect to the product decomposition, i.e.,
$T \in \Gamma(\Lambda^2 V_1^* \otimes V_2^*) \quad \text{or} \quad T \in \Gamma(V_1^* \otimes \Lambda^2 V_2^*)$ and assume further that the harmonic representative $[\omega_h]$ is a non‑trivial mixed class in the Künneth decomposition of $H^3(M)$.

Define the mixed factor subspaces
\[
\mathcal{V}_{2,1}:=\operatorname{span}\{\alpha_i\otimes\beta_i\}\subset\mathcal{H}^2(M_1)\otimes\mathcal{H}^1(M_2),
\]
\[
\mathcal{V}_{1,2}:=\operatorname{span}\{\tilde\alpha_j\otimes\tilde\beta_j\}\subset\mathcal{H}^1(M_1)\otimes\mathcal{H}^2(M_2),
\]
where the spans are taken over any minimal harmonic decomposition of the mixed Künneth 
components of $[\omega]$; in the admissible setting these spans are independent of the 
particular choice (Lemma \ref{lem:invariance}).

Define the obstruction kernels
\[
\ker\Psi_p:=\mathcal{V}_{2,1}\cap\bigl(\mathcal{P}_2(M_1)\otimes\mathcal{P}_1(M_2)\bigr),\quad
\ker\widetilde{\Psi}_p:=\mathcal{V}_{1,2}\cap\bigl(\mathcal{P}_1(M_1)\otimes\mathcal{P}_2(M_2)\bigr),
\]
and set $\mathcal{K}:=\ker\Psi_p\oplus\ker\widetilde{\Psi}_p$.

Then, there exists a non-empty open subset $\mathcal{V} \subset M$ such that for every $p\in \mathcal{V}$,
\[
\dim\bigl(\mathfrak{hol}_p^{\mathrm{off}}(\nabla^{C})\bigr)\;\geq\;
r^\sharp\;:=\;\operatorname{rank}_{\mathbb{R}}\bigl([\omega]_{\mathrm{mixed}}\bigr)-\dim\mathcal{K}.
\]
The bound depends only on $[\omega]_{\text{mixed}}$ and on the parallel-form strata $\mathcal{P}_k(M_i)$; hence it is invariant under any deformation of $g_1$ or $g_2$ that preserves these strata.
\end{theorem}

\begin{proof}
Let $\omega_h$ be the $g$-harmonic representative of $[\omega]$. Decompose $\omega_h$ according to K\"unneth and Hodge on the factors. For the $(2,1)$-part, write
\begin{equation}\label{eq:decomp21}
\omega_h^{2,1}=\sum_{i=1}^{r_{2,1}} \alpha_i\wedge \beta_i,
\end{equation}
with $\alpha_i\in \mathcal{H}^2(M_1)$ and $\beta_i\in \mathcal{H}^1(M_2)$ harmonic and chosen so that the number of terms is minimal. Similarly,
\begin{equation}\label{eq:decomp12}
\omega_h^{1,2}=\sum_{j=1}^{r_{1,2}} \tilde\alpha_j\wedge \tilde\beta_j,
\end{equation}
with $\tilde\alpha_j\in \mathcal{H}^1(M_1)$ and $\tilde\beta_j\in \mathcal{H}^2(M_2)$. The total number of simple mixed tensors equals $r:=r_{2,1}+r_{1,2}=\operatorname{rank}_\mathbb{R}\big([\omega]_{\mathrm{mixed}}\big)$. Under the pure bidegree hypothesis, exactly one of $r_{2,1}$, $r_{1,2}$ is non-zero; the other vanishes together with the corresponding decomposition.

Set $r^\sharp\;:=\;\operatorname{rank}_{\mathbb{R}}\bigl([\omega]_{\mathrm{mixed}}\bigr)-\dim\mathcal{K}$. We will prove $\dim\mathfrak{hol}_{p}^{\operatorname{off}}(\nabla^{C})\geq r^\sharp$.

Consider the curvature formula \eqref{eq:curv-torsion}
for $X\in V_1$ and $Y\in V_2$:
\[
  R^C(X,Y)
  \;=\;
  R^{LC}(X,Y)
  \;+\;
  \tfrac{1}{2}\bigl[
    (\nabla^{LC}_X T)(Y,\cdot\,,\cdot)^{\sharp}
    - (\nabla^{LC}_Y T)(X,\cdot\,,\cdot)^{\sharp}
  \bigr]
  \;+\;
  \tfrac{1}{4}\,Q_T(X,Y).
\]
On a Riemannian product the mixed Riemann curvature
vanishes: $R^{LC}(X,Y)=0$ whenever $X\in V_1$
and $Y\in V_2$
(since $\nabla^{LC}_XY=0=[X,Y]$ for vector fields
tangent to different factors).
For the remaining three terms we exploit the rigid
pointwise algebraic structure that the surface
hypothesis imposes on mixed $3$-forms.

Since both factors are compact oriented surfaces,
$\dim V_i=2$, and the exterior power
$\Lambda^2V_1^*$ is one-dimensional, spanned by the
area form $\operatorname{vol}_{V_1}$.
Every $3$-form of bidegree $(2,1)$ on $M$ can
therefore be written at each point $p$ as
\begin{equation}\label{eq:T21-ptwise}
  T^{2,1}_p
  \;=\;
  h(p)\;\operatorname{vol}_{V_1}\wedge\tau_p,
\end{equation}
where $h\colon M\to\mathbb{R}$ is a smooth function
and $\tau$ is a smooth section of $V_2^*$
(depending on the full point $p=(p_1,p_2)$).
Likewise, every $(1,2)$-form satisfies
$T^{1,2}_p
= \widetilde{h}(p)\;\widetilde{\sigma}_p
\wedge\operatorname{vol}_{V_2}$
with $\widetilde{\sigma}\in\Gamma(V_1^*)$.
These representations hold for the \emph{full}
torsion $T$, not merely for its harmonic
projection $\omega_h$.

We denote by
$J_i\colon V_i\to V_i$ the volume endomorphism of
the oriented surface $M_i$, characterised by
$g(J_iu,v)=\operatorname{vol}_{V_i}(u,v)$;
it is a skew-symmetric isometry with
$J_i^2=-\operatorname{Id}$, and in particular
$J_iX\perp X$ for every non-zero $X\in V_i$.
Let
$\pi_0\colon\mathfrak{so}(T_pM)
\to\mathfrak{g}_0
:=\mathfrak{so}(V_1)\oplus\mathfrak{so}(V_2)$
and
$\pi_1\colon\mathfrak{so}(T_pM)
\to\mathfrak{g}_1
:=\mathfrak{so}^{\mathrm{off}}(T_pM)$
denote the projections induced by the
$\mathbb{Z}_2$-grading of $\mathfrak{so}(T_pM)$.

A direct computation
using \eqref{eq:T21-ptwise}
yields the following structure for the three
surviving terms.
The first,
$\frac{1}{2}(\nabla^{LC}_X T)(Y,\cdot\,,\cdot)^{\sharp}$,
is purely diagonal: applying the Leibniz rule to
$\nabla^{LC}_X(h\,\operatorname{vol}_{V_1}\wedge\tau)$
and using
$\nabla^{LC}_X\operatorname{vol}_{V_1}=0$
(since $X\in V_1$ and the volume form of $M_1$ is
parallel), one finds that this term is proportional
to the volume endomorphism $J_1$ on $V_1$
(respectively $J_2$ on $V_2$ for the
$(1,2)$-component of $T$); thus
\begin{equation}\label{eq:diag-comp}
  \pi_0\bigl(R^C(X,Y)\bigr)
  \;=\;
  \lambda(X,Y)\,J_1 \;+\; \mu(X,Y)\,J_2
\end{equation}
for certain scalars $\lambda,\mu$. The scalar $\lambda$ arises from the $(2,1)$-component of $T$ and its covariant derivatives, while $\mu$ arises from the $(1,2)$-component; under the pure bidegree hypothesis exactly one of $\lambda$, $\mu$ is identically zero.
The second,
$-\frac{1}{2}(\nabla^{LC}_YT)(X,\cdot\,,\cdot)^{\sharp}$,
is purely off-diagonal.
Defining the $1$-form
$\rho_Y\in V_2^*$ by
\begin{equation}\label{eq:rho-def}
  \rho_Y
  \;:=\;
  (Yh)\,\tau
  \;+\;
  h\,\nabla^{g_2}_Y\!\tau
  \;=\;
  \nabla^{g_2}_Y(h\tau),
\end{equation}
the restriction to
$\operatorname{Hom}(V_1,V_2)$ reads
\begin{equation}\label{eq:B-term}
  -\tfrac{1}{2}
  (\nabla^{LC}_YT^{2,1})(X,Z,\cdot)^{\sharp}
  \Big|_{Z\in V_1}
  \;=\;
  -\tfrac{1}{2}\,
  g(J_1X,Z)\;\rho_Y^{\sharp}
  \;\in V_2.
\end{equation}
The third, $\frac{1}{4}Q_T(X,Y)$, is also
purely off-diagonal (by a bidegree computation),
with restriction
\begin{equation}\label{eq:Q-term}
  \tfrac{1}{4}\,Q_T(X,Y)\big|_{Z\in V_1}
  \;=\;
  \tfrac{h^2}{4}\;\tau(Y)\;
  g(X,Z)\;\tau^{\sharp}
  \;\in V_2.
\end{equation}

To derive \eqref{eq:Q-term}, observe that the contortion tensor $K_XY = \tfrac{1}{2}T(X,Y,\cdot)^\sharp$ acts on $Z \in V_1$ via the $(2,1)$-component as
\[
K_Y^{(2,1)}Z = \frac{h}{2}\,\tau(Y)\,J_1 Z,
\]
since $(\operatorname{vol}_{V_1} \wedge \tau)(Y,Z,W) = \tau(Y)\,\operatorname{vol}_{V_1}(Z,W)$ for $Y \in V_2$ and $Z,W \in V_1$. The iterated contortion $K_X(K_YZ)$ with $X \in V_1$ then gives
\[
K_X(K_YZ) = \frac{h^2}{4}\,\tau(Y)\,g(X,Z)\,\tau^\sharp,
\]
where we used $T^{(2,1)}(X,J_1Z,W) = h\,\operatorname{vol}_{V_1}(X,J_1Z)\,\tau(W) = h\,g(J_1X,J_1Z)\,\tau(W) = h\,g(X,Z)\,\tau(W)$, the last equality following from $J_1$ being an isometry. Conversely, $K_Y(K_XZ) = 0$: since $K_XZ = \tfrac{h}{2}\operatorname{vol}_{V_1}(X,Z)\,\tau^\sharp \in V_2$ and $T^{(2,1)}(Y,W,\cdot) = 0$ whenever $Y,W \in V_2$ (the form $\operatorname{vol}_{V_1}$ requires two arguments in $V_1$). Therefore $\tfrac14 Q_T(X,Y)Z = K_X(K_YZ) - K_Y(K_XZ)$ yields \eqref{eq:Q-term}.

Combining \eqref{eq:B-term}
and \eqref{eq:Q-term}, the full off-diagonal
projection restricted to
$\operatorname{Hom}(V_1,V_2)$ is
\begin{equation}\label{eq:pi1-full}
  \pi_1\bigl(R^C(X,Y)\bigr)\big|_{V_1\to V_2}
  \;:\;
  Z\;\longmapsto\;
  -\tfrac{1}{2}\,g(J_1X,Z)\;\rho_Y^{\sharp}
  \;+\;
  \tfrac{h^2}{4}\,\tau(Y)\,g(X,Z)\;\tau^{\sharp}.
\end{equation}
The two summands act through the functionals
$g(J_1X,\cdot\,)$ and $g(X,\cdot\,)$ on $V_1$
respectively.
Because $J_1X\perp X$, these functionals are
orthogonal, and the two contributions cannot
cancel: if the map \eqref{eq:pi1-full}
vanishes for all $Z\in V_1$, then
\begin{equation}\label{eq:no-cancel}
  \rho_Y^{\sharp}=0
  \qquad\text{and}\qquad
  h^2\,\tau(Y)\,\tau^{\sharp}=0
\end{equation}
must hold separately.

We now observe that the quadratic
term \eqref{eq:Q-term} is non-vanishing on a
non-empty open set.
If $r_{2,1}=1$, the harmonic projection
$\omega_h^{2,1}
=c\cdot\operatorname{vol}_{M_1}\wedge\beta$
is non-zero, so
$[T^{2,1}]=[\omega_h^{2,1}]\neq 0$
in cohomology.
Since $T^{2,1}=\omega_h^{2,1}+\eta^{2,1}$
with $\eta^{2,1}$ exact plus coexact,
Hodge orthogonality gives
$T^{2,1}\not\equiv 0$ on $M$
(otherwise $\omega_h^{2,1}=-\eta^{2,1}$ would
be simultaneously harmonic and
exact-plus-coexact, forcing
$\omega_h^{2,1}=0$).
The set
\begin{equation}\label{eq:W21-def}
  \mathcal{W}^{2,1}
  \;:=\;
  \{p\in M:T^{2,1}_p\neq 0\}
\end{equation}
is therefore a non-empty open subset of $M$.
On $\mathcal{W}^{2,1}$ we have
$h(p)\neq 0$ and $\tau_p\neq 0$
(by \eqref{eq:T21-ptwise}), so choosing
$Y\in V_2$ with $\tau_p(Y)\neq 0$
yields
$h^2\tau(Y)\tau^{\sharp}\neq 0$;
the second condition in \eqref{eq:no-cancel}
is therefore violated, and
$\pi_1(R^C(X,Y))\neq 0$.
When $T$ has pure bidegree $(1,2)$, the symmetric argument for the $(1,2)$-component (using $J_2Y\perp Y$) produces a non-empty open set $\mathcal{W}^{1,2}$ whenever $r_{1,2}=1$. Under the pure bidegree hypothesis, at most one of the sets $\mathcal{W}^{2,1}$, $\mathcal{W}^{1,2}$ is non-empty; both being simultaneously non-empty would require $T$ to have mixed bidegree, contradicting the hypothesis.

In particular, the orthogonality of the two summands 
in \eqref{eq:pi1-full} ensures that, on $\mathcal{W}^{2,1}$, the off-diagonal projection $\pi_1(R^C(X,Y))$ is non-zero regardless of the exact-plus-coexact part $\eta = T - \omega_h$. Indeed, $h$ and $\tau$ in \eqref{eq:T21-ptwise} are 
determined by the full torsion $T$, and on $\mathcal{W}^{2,1}$ one has $h(p)\neq 0$ and $\tau_p\neq 0$ by definition; the quadratic term $\tfrac{h^2}{4}\tau(Y)\,\tau^\sharp$ is therefore non-zero for a suitable choice of $Y\in V_2$, and acts through the functional $g(X,\cdot\,)$, which is orthogonal to $g(J_1X,\cdot\,)$ appearing in the linear term.  Thus, even if $\eta$ alters the pointwise values of $\rho_Y$, $h$, and $\tau$ relative to those of the harmonic representative $\omega_h$, the two contributions to the off-diagonal projection cannot cancel.  In the admissible setting, 
$r^\sharp\leq 1$ per bidegree component, so this single protected direction suffices: the cohomological count $r^\sharp$ records these directions; the purely off-diagonal holonomy element itself is produced in \emph{Case 1}. The explicit construction of a \emph{purely} off-diagonal holonomy element from the full 
curvature $R^C$ is carried out in \emph{Case $1$} and \emph{Case $2$} below.

Under the pure bidegree hypothesis exactly one of $r_{2,1}$, $r_{1,2}$ is non-zero; we treat the $(2,1)$ component, the $(1,2)$ case being symmetric. To count the mixed directions not absorbed by the parallel-form strata, we introduce a linear map and compute its kernel.

Define the linear map
\begin{equation}\label{eq:defPsi}
\Psi_p:\mathcal{V}_{2,1}\longrightarrow T_p^*M\otimes\Lambda^3 T_p^*M,\qquad
\Psi_p\Bigl(\sum_{i=1}^{r_{2,1}} c_i\,\alpha_i\otimes\beta_i\Bigr):=
\sum_{i=1}^{r_{2,1}} c_i\,\bigl(\nabla^{LC}(\alpha_i\wedge\beta_i)\bigr)\big|_p .
\end{equation}
Here $\nabla^{LC}(\alpha_i\wedge\beta_i)$ is the covariant derivative of the harmonic $3$-form $\alpha_i\wedge\beta_i$, a section of $T^*M\otimes\Lambda^3 T^*M$, evaluated at $p$. Thus $\Psi_p$ is a genuine linear map between finite-dimensional
vector spaces, defined without any auxiliary choice of vectors; its sole purpose is to compute $\dim\ker\Psi_p$, hence $\dim\mathcal{K}$. The elements of
$\operatorname{im}\Psi_p$ are covariant derivatives of harmonic forms, \emph{not} curvature operators; the purely off-diagonal holonomy element is produced separately, in \emph{Case 1} below, from the full curvature $R^C(X_0,Y)$, which is
a genuine skew-symmetric endomorphism of $T_pM$.

\medskip
\noindent\textbf{Claim 5.2.1 (Kernel of $\Psi_p$).}
On a non-empty open subset of $M$ one has
\[
\ker\Psi_p=\mathcal{V}_{2,1}\cap\bigl(\mathcal{P}_2(M_1)\otimes\mathcal{P}_1(M_2)\bigr).
\]

\begin{proof}[Proof of Claim]
If $r_{2,1}=0$ both sides are $\{0\}$. Assume $r_{2,1}=1$ and write
$\xi=c_1(\alpha_1\otimes\beta_1)$. By the Leibniz rule on the Riemannian product,
\[
\nabla^{LC}(\alpha_1\wedge\beta_1)
=(\nabla^{g_1}\alpha_1)\wedge\beta_1+\alpha_1\wedge(\nabla^{g_2}\beta_1).
\]
Since $M_1$ is a compact oriented surface,
$\alpha_1\in\mathcal{H}^2(M_1)=\mathbb{R}\cdot\operatorname{vol}_{M_1}$ is parallel,
so $\nabla^{g_1}\alpha_1=0$ and
\[
\Psi_p(\xi)=c_1\,\alpha_1\wedge(\nabla^{g_2}\beta_1)\big|_{p}.
\]
Since $\operatorname{vol}_{M_1}$ is a nowhere-vanishing section of $\Lambda^2 V_1^*$,
the $3$-form $\alpha_1\wedge(\nabla^{g_2}_W\beta_1)$ vanishes for every $W\in V_2$ if
and only if $\nabla^{g_2}\beta_1=0$ at $p_2$. Hence $\Psi_p(\xi)=0$ if and only if
$c_1=0$ or $(\nabla^{g_2}\beta_1)(p_2)=0$. Two cases arise.

If $\beta_1\in\mathcal{P}_1(M_2)$, then $\nabla^{g_2}\beta_1\equiv 0$, so $\Psi_p\equiv 0$
and $\ker\Psi_p=\mathcal{V}_{2,1}=\operatorname{span}\{\alpha_1\otimes\beta_1\}$,
which coincides with
$\mathcal{V}_{2,1}\cap(\mathcal{P}_2(M_1)\otimes\mathcal{P}_1(M_2))$ because
$\alpha_1\in\mathcal{P}_2(M_1)$ and $\beta_1\in\mathcal{P}_1(M_2)$. The identity holds
on all of $M$.

If $\beta_1\notin\mathcal{P}_1(M_2)$, then $\nabla^{g_2}\beta_1\not\equiv 0$, and by
Lemma \ref{lem:maximal-rank-open} --- which, as noted after its proof, holds for
harmonic forms of any degree on either factor, hence for the harmonic $1$-form
$\beta_1$ on $M_2$ --- the set
$\{p_2\in M_2:(\nabla^{g_2}\beta_1)(p_2)\neq 0\}$ is a non-empty open subset of
$M_2$. For every $p=(p_1,p_2)$ with $p_2$ in this set, $\Psi_p(\xi)=0$ if and only if
$c_1=0$, so $\ker\Psi_p=\{0\}$, which coincides with
$\mathcal{V}_{2,1}\cap(\mathcal{P}_2(M_1)\otimes\mathcal{P}_1(M_2))=\{0\}$ (as
$\beta_1\notin\mathcal{P}_1(M_2)$). The identity holds on the non-empty open set
$M_1\times\{p_2:(\nabla^{g_2}\beta_1)(p_2)\neq 0\}$.
\end{proof}

Denote by $\mathcal{O}\subseteq M$ the non-empty open set provided by the Claim (all
of $M$ in the parallel case; $M_1\times\{(\nabla^{g_2}\beta_1)\neq 0\}$ otherwise).
For $p\in\mathcal{O}$ the induced map on the quotient,
\[
\overline{\Psi}_p:\mathcal{V}_{2,1}/\ker\Psi_p\longrightarrow
T_p^*M\otimes\Lambda^3 T_p^*M,\qquad \overline{\Psi}_p([\xi]):=\Psi_p(\xi),
\]
is well-defined, linear and injective; hence
\begin{equation}\label{eq:dimim}
\dim(\operatorname{im}\Psi_p)=\dim\bigl(\mathcal{V}_{2,1}/\ker\Psi_p\bigr)
=r_{2,1}-\dim\ker\Psi_p .
\end{equation}
An identical argument for the $(1,2)$ component gives, on a corresponding non-empty
open subset,
\begin{equation}\label{eq:dimim12}
\dim(\operatorname{im}\widetilde{\Psi}_p)=r_{1,2}-\dim\ker\widetilde{\Psi}_p .
\end{equation}
Consequently, on a non-empty open subset of $M$,
\[
\dim(\operatorname{im}\Psi_p)+\dim(\operatorname{im}\widetilde{\Psi}_p)
=\bigl(r_{2,1}-\dim\ker\Psi_p\bigr)+\bigl(r_{1,2}-\dim\ker\widetilde{\Psi}_p\bigr)
=r-\dim\mathcal{K}=r^\sharp,
\]
with exactly one summand non-zero under the pure bidegree hypothesis. This is a purely cohomological count of the mixed simple-tensor directions not absorbed by the parallel-form strata; the elements of $\operatorname{im}\Psi_p$,
$\operatorname{im}\widetilde{\Psi}_p$ are covariant derivatives of harmonic $3$-forms, not curvature operators.

It remains to produce a purely off-diagonal element of the holonomy algebra.  Under the pure bidegree hypothesis exactly one of $r_{2,1}$, $r_{1,2}$ is non-zero, and for admissible manifolds each is at most $1$; consequently $r^{\sharp}=r-\dim\mathcal{K}\leq 1$, and it suffices to exhibit a single purely off-diagonal holonomy element on a non-empty open set.

Without loss of generality assume $T$ has pure bidegree $(2,1)$, so that $r_{2,1} \geq 1$ and $r_{1,2} = 0$ (the case of pure bidegree $(1,2)$ is entirely symmetric).  We work on the non-empty open set $\mathcal{V} := \mathcal{W}^{2,1}$, where $\mathcal{W}^{2,1}$ is the
non-empty open set introduced in \ref{eq:W21-def}. At every point $p \in \mathcal{V}$ the torsion satisfies $h(p) \neq 0$ and $\tau_p \neq 0$ by definition of $\mathcal{W}^{2,1}$.

Consider the mixed curvature span
\[
  \mathcal{S}_p
  \;:=\;
  \operatorname{span}\bigl\{
    R^C(X,Y) : X\in V_1,\; Y\in V_2
  \bigr\}.
\]
By the Ambrose--Singer theorem \cite{A-S (1953)},
$\mathcal{S}_p
\subseteq\mathfrak{hol}_p(\nabla^C)$.

We now show that $\mathfrak{hol}_p(\nabla^C)$ contains
at least one purely off-diagonal element on a non-empty open
subset of $M$.

Suppose $r_{2,1}=1$.
Recall that at every point $p\in\mathcal{W}^{2,1}$ we have
$h(p)\neq 0$ and $\tau_p\neq 0$.
We distinguish two cases according to whether $p$ lies in the
interior of $\{T^{1,2}=0\}$.

\medskip
\noindent\textit{Case 1: $p\in\mathcal{W}^{2,1}\setminus
\overline{\mathcal{W}^{1,2}}$.}
Since $\mathcal{W}^{1,2}=\{q\in M:T^{1,2}_q\neq 0\}$ is open,
its complement $\{T^{1,2}=0\}$ is closed and its interior is
$M\setminus\overline{\mathcal{W}^{1,2}}$.  At every point
$p\in\mathcal{W}^{2,1}\setminus\overline{\mathcal{W}^{1,2}}$
the $(1,2)$-component $T^{1,2}$ vanishes identically on an open
neighbourhood of $p$; consequently
$\nabla^{\mathrm{LC}}T^{1,2}\big|_p=0$ as well
(the derivative of a section that vanishes on an open set is zero
on that set).

It follows that both the quadratic torsion terms involving $T^{1,2}$
and the $\nabla^{\mathrm{LC}}T^{1,2}$-terms in the curvature
formula \eqref{eq:curv-torsion} vanish at $p$.
The diagonal projection therefore reduces to
\[
  \pi_0\bigl(R^C(X,Y)\bigr) = \lambda(X,Y)\,J_1
\]
(the $\mu(X,Y)\,J_2$ contribution is absent because it arises
entirely from $T^{1,2}$ and its covariant derivatives, both of which
vanish at $p$).
Choose $Y\in V_2$ with $\tau_p(Y)\neq 0$ (possible because
$p\in\mathcal{W}^{2,1}$ implies $\tau_p\neq 0$).  The linear
functional $\lambda(\cdot,Y)\colon V_1\to\mathbb{R}$ has kernel of
dimension at least $1$ (since $\dim V_1=2$).  Pick
$X_0\in\ker\bigl(\lambda(\cdot,Y)\bigr)\setminus\{0\}$.  Then:
\begin{itemize}
\item [] \emph{Diagonal part:}
  $\pi_0\bigl(R^C(X_0,Y)\bigr)=\lambda(X_0,Y)\,J_1=0$.
\medskip
\item [] \emph{Off-diagonal part evaluated on $Z=X_0$:}
  by \eqref{eq:pi1-full}, the two contributions are
  $-\tfrac{1}{2}\,g(J_1X_0,X_0)\,\rho_Y^\sharp$ and
  $\tfrac{h^2}{4}\,\tau(Y)\,g(X_0,X_0)\,\tau^\sharp$.
  The first vanishes because $g(J_1X_0,X_0)=0$
  (skew-symmetry of $J_1$); the second equals
  $\tfrac{h^2}{4}\,\tau(Y)\,|X_0|^2\,\tau^\sharp\neq 0$
  (since $h\neq 0$, $\tau(Y)\neq 0$, and $|X_0|^2>0$).
\end{itemize}
Therefore $R^C(X_0, Y)$ is a non-zero element of 
$\mathfrak{so}^{\mathrm{off}}(T_pM)$.  By the Ambrose--Singer 
theorem \cite{A-S (1953)}, $R^C(X_0, Y) \in \mathfrak{hol}_p(\nabla^C)$.  
Since $\mathfrak{hol}^{\mathrm{off}}_p(\nabla^C) 
= \mathfrak{hol}_p(\nabla^C) \cap 
\mathfrak{so}^{\mathrm{off}}(T_pM)$ (Definition \ref{def:off-hol}), we conclude that $R^C(X_0, Y) \in 
\mathfrak{hol}^{\mathrm{off}}_p(\nabla^C)$, whence 
$\dim\mathfrak{hol}^{\mathrm{off}}_p(\nabla^C) \geq 1$.

\medskip

\noindent \textit{Case 2: $p \in \mathcal{W}^{2,1} \cap \mathcal{W}^{1,2}$ (both components present).} This case is vacuous under the pure bidegree hypothesis. We show that the hypothesis
$T\in\Gamma(\Lambda^{2}V_{1}^{*}\otimes V_{2}^{*})$
(the argument for $T\in\Gamma(V_{1}^{*}\otimes\Lambda^{2}V_{2}^{*})$ is
entirely symmetric) eliminates the overlap
$\mathcal{W}^{2,1}\cap\mathcal{W}^{1,2}$ altogether.

The pure bidegree condition requires that for every $p\in M$ and all
$X\in V_{1}(p)$, $Y_{1},Y_{2}\in V_{2}(p)$ one has
$T_{p}(X,Y_{1},Y_{2})=0$,
which is precisely the statement $T^{1,2}\equiv 0$ as a section of
$V_{1}^{*}\otimes\Lambda^{2}V_{2}^{*}$.  Consequently
\[
  \mathcal{W}^{1,2}
  \;=\;\bigl\{p\in M : T^{1,2}(p)\neq 0\bigr\}
  \;=\;\emptyset,
\]
and therefore $\mathcal{W}^{2,1}\cap\mathcal{W}^{1,2}=\emptyset$.
This is not a generic or measure-theoretic statement: the overlap is
empty at every point of $M$. In particular, every
$p\in\mathcal{W}^{2,1}$ admits a neighbourhood on which $T^{1,2}$ vanishes identically, so the hypotheses of Case 1 are satisfied at all such points.

Conversely, if $\mathcal{W}^{2,1}\cap\mathcal{W}^{1,2}\neq\emptyset$
then there exists a point at which both $T^{2,1}$ and $T^{1,2}$ are non-zero, whence $T$ cannot have pure bidegree. The pure bidegree hypothesis is therefore the precise condition that renders Case 2
vacuous.

As established in \eqref{eq:W21-def} above via Hodge orthogonality, $\mathcal{W}^{2,1}$ is a non-empty open subset of $M$.

Recall $\mathcal{V}:=\mathcal{W}^{2,1}$ (in the symmetric case $T = T^{1,2}$ one sets $\mathcal{V} := \mathcal{W}^{1,2}$).  This is a non-empty open subset
of $M$, and at every point of $\mathcal{V}$ the hypotheses of \emph{Case 1} are in force (since $\mathcal{W}^{1,2}=\emptyset$).  Applying the
conclusion of Case 1 we obtain, for every $p\in\mathcal{V}$,
\[
\mathfrak{hol}_{p}^{\mathrm{off}}(\nabla^{C})\neq\{0\},
\]
together with the rank estimate
\[
\dim\bigl(\mathfrak{hol}_p^{\mathrm{off}}(\nabla^{C})\bigr)\;\geq\;
r^\sharp\;:=\;\operatorname{rank}_{\mathbb{R}}\bigl([\omega]_{\mathrm{mixed}}\bigr)-\dim\mathcal{K}.
\]
The quantity $r^{\sharp}$ depends only on the mixed K\"unneth class
$[\omega]_{\mathrm{mixed}}$ and on the parallel-form strata
$\mathcal{P}_{k}(M_{i})$; it is therefore invariant under any
deformation of $g_{1}$ or $g_{2}$ that preserves these strata.
    
\end{proof}

\begin{corollary} Under the hypotheses of Theorem \ref{thm:main}, if $\operatorname{rank}_\mathbb{R}([\omega]_{\mathrm{mixed}}) \geq 1$, then $\dim\mathfrak{hol}^{\mathrm{off}}_p(\nabla^C) \geq 1$ on the non-empty open subset $\mathcal{V} \subset M$ of Theorem \ref{thm:main}. In particular, the product splitting $V_1 \oplus V_2$ is not invariant under the holonomy of $\nabla^C$; equivalently, $\nabla^C$ does not reduce to a product connection along $V_1\oplus V_2$.
\end{corollary}

\begin{proof} Immediate from \emph{Case 1} of the proof of Theorem \ref{thm:main}: the quadratic torsion term $\frac{1}{4}Q_T$ produces a non-zero purely off-diagonal element of $\mathfrak{hol}_p(\nabla^C)$ at every point of $\mathcal{W}^{2,1} \neq \emptyset$ (respectively $\mathcal{W}^{1,2} \neq \emptyset$), using only the non-vanishing of $T^{2,1}$ (respectively $T^{1,2}$), which is guaranteed by $\operatorname{rank}_\mathbb{R}([\omega]_{\mathrm{mixed}}) \geq 1$.
\end{proof}

\begin{remark}[Role of the harmonic projection]
\label{rem:harmonic}
The invariant $r^\sharp$ is defined through the harmonic projection $\omega_h$ alone, not through the full torsion 
$T = \omega_h + d\alpha + \delta\beta$.  This is essential for invariance: for a fixed class $[\omega]$, infinitely many torsion forms represent $[\omega]$ (corresponding to 
different exact and coexact parts), so any quantity depending on the full torsion would not be an invariant of the cohomology class.

As a consequence, $r^\sharp$ captures only the off-diagonal holonomy forced by the cohomological datum.
More precisely, the linear map $\Psi_p$ defined in \eqref{eq:defPsi} serves only to \emph{count} the number of mixed simple-tensor directions that are not parallel: its
injectivity on the quotient $\mathcal{V}_{2,1}/\ker\Psi_p$, which holds on a non-empty open subset of $M$ (see the proof of Theorem \ref{thm:main}), yields
$\dim\operatorname{im}\Psi_p = r_{2,1}-\dim\ker\Psi_p$. The elements of $\operatorname{im}\Psi_p$ are covariant derivatives of harmonic $3$-forms, not curvature operators.
The \emph{existence} of a purely off-diagonal holonomy element is then established by \emph{Case 1} of the proof through the quadratic torsion term $\tfrac{1}{4}Q_T$, which depends on the pointwise values of $T$ rather than on its cohomology class alone.
Since $Q_T$ may generate off-diagonal curvature beyond what is detected by $r^\sharp$, the invariant $r^\sharp$ provides a lower bound, not an exact count, for the off-diagonal holonomy dimension. More precisely, $r^\sharp$ certifies the non-preservation of the product splitting $V_1\oplus V_2$ that arises when the de-Rham--irreducible geometry of a factor prevents the cohomological coupling from being absorbed by the parallel-form strata; when $r^\sharp = 0$, no such obstruction to absorption is present, and any failure of $V_1\oplus V_2$ to be holonomy-invariant, though it may still occur via the quadratic torsion term $Q_T$ (as discussed above), does not arise from this mechanism (see Remark \ref{rem:deRham}). In this sense, $r^\sharp$ quantifies the \emph{transverse non-reducibility of the connection} (Remark \ref{rem:deRham}).

However, in the generic situation where $\nabla^{LC}T \neq 0$, the curvature operator $R^C(X,Y)$ has both diagonal and off-diagonal components, and the off-diagonal projection of a holonomy element is not itself a holonomy element in general; extracting a \emph{purely} off-diagonal element from the full curvature requires the algebraic mechanism of \emph{Theorem \ref{thm:main}}.
\end{remark}

\begin{remark} [$r^\sharp$ as a property of the pair $(M,g)$]
\label{rem:property}
The condition $r^\sharp([\omega], g) \geq 1$ for some class $[\omega] \in H^3(M)$ is a property of the pair $(M,g)$ in admissible manifolds. When satisfied, $r^\sharp$ provides a non-trivial lower bound $\dim\mathfrak{hol}^{\mathrm{off}}_p \geq 1$, certifying that, for every metric connection with totally antisymmetric torsion of pure bidegree representing $[\omega]$, the product splitting $V_1\oplus V_2$ fails to be holonomy-invariant, an obstruction that cannot be removed by metric deformations preserving the parallel-form strata. When $r^\sharp([\omega], g) = 0$ for every class $[\omega]$, no such assurance is provided: for instance, on $S^2 \times S^2$ the group $H^3(M) = 0$, so no metric connection with totally antisymmetric torsion can represent a non-trivial cohomology class.
\end{remark}

\begin{remark}[Geometric interpretation of $r^\sharp$ and the obstruction space $\mathcal{K}$]
\label{rem:deRham}
The invariant $r^{\sharp} = \operatorname{rank}_{\mathbb{R}}\bigl([\omega]_{\mathrm{mixed}}\bigr) - \dim\mathcal{K}$ is computable from the cohomology class $[\omega_h]$ and the parallel-form strata $\mathcal{P}_k(M_i)$ alone, without knowledge of the specific torsion representative $T$.
The correction $-\dim\mathcal{K}$ admits a natural geometric
interpretation in terms of the de Rham decomposition of the factors.

When $\dim\mathcal{K}>0$, the subspace $\mathcal{K}$ contains at least one simple tensor $\alpha_k\otimes\beta_k$ with
$\alpha_k\in\mathcal{P}_2(M_i)$ and $\beta_k\in\mathcal{P}_1(M_j)$, i.e. both factors are parallel with respect to the Levi-Civita connections of the factor metrics.
In the admissible setting, the condition on $\alpha_k$ is automatic ($\mathcal{P}_2(M_i)=\mathbb{R}\cdot\operatorname{vol}_{M_i}$ for every compact oriented surface), and the operative restriction concerns only $\beta_k$.
A non-zero parallel $1$-form $\beta_k$ on a compact Riemannian
manifold defines a parallel vector field $\beta_k^{\sharp}$, which by the de Rham decomposition theorem determines a flat sub-factor in the universal cover: $\widetilde{M}_j$ splits isometrically as $\widetilde{M}_j\cong\widetilde{N}_j\times\mathbb{R}_{\beta_k}$.
In the admissible setting, $\mathcal{P}_1(M_j)\neq\{0\}$ forces
$M_j$ to be a flat torus, and the splitting descends to the
manifold itself: $M_j\cong N_j\times S^1_{\beta_k}$.
The component $\alpha_k\otimes\beta_k$ of the torsion class then couples $M_i$ exclusively with the flat direction $S^1_{\beta_k}$, without involving $N_j$.

The components not in $\mathcal{K}$ (those with
$\beta_l\notin\mathcal{P}_1(M_j)$, or
$\alpha_l\notin\mathcal{P}_2(M_i)$ in higher-dimensional settings) couple $M_i$ with the de Rham irreducible part of $M_j$, the geometrically non-trivial part of the factor.

The invariant $r^{\sharp}$ therefore measures the dimension of the part of the mixed cohomological coupling that is \emph{transverse} to the de Rham splitting, i.e. the part that involves the irreducible factor.
The correction $\dim\mathcal{K}$ measures the dimension of the part \emph{aligned} with the flat sub-factors.

This interpretation clarifies the meaning of the two regimes:
\begin{itemize}
\item[(i)] When $r^{\sharp}\geq 1$: the irreducible de Rham geometry of a factor prevents the cohomological coupling from being absorbed by the parallel-form strata. Theorem \ref{thm:main} certifies that this non-absorption forces the product splitting $V_1\oplus V_2$ to fail to be holonomy-invariant, and this failure cannot be removed by metric deformations preserving the parallel-form strata. We refer to this phenomenon as the \emph{transverse non-reducibility of the connection}: the holonomy fails to preserve $V_1\oplus V_2$ precisely because of the de-Rham--irreducible geometry of a factor: see Proposition \ref{prop:reduction} for a setting in which this transversality is detected by $\dim\mathcal{K}$ alone.
\item[(ii)] When $r^{\sharp}=0$: no irreducible de Rham factor prevents absorption; the entire cohomological coupling is absorbed by the parallel-form strata. Transverse non-reducibility is absent. This does \emph{not} mean that $V_1\oplus V_2$ is holonomy-invariant: the quadratic torsion term $\frac{1}{4}Q_T$ can still generate off-diagonal curvature (Remark \ref{rem:harmonic}), so the splitting may still fail to be preserved; it only signals that this failure does not arise from the mechanism detected by $r^\sharp$.
\end{itemize}

In the admissible setting (products of surfaces), the distinction reduces to the dichotomy of Remark \ref{rmk:dichotomy}.
For higher-dimensional factors, the de Rham decomposition is richer: $k$ independent parallel $1$-forms on $M_j$ induce a splitting of the universal cover $\widetilde{M}_j\cong\widetilde{N}_j\times\mathbb{R}^k$;
when this splitting descends to the manifold (e.g.\ under a simple connectivity hypothesis or for factors with special holonomy), one obtains $M_j\cong N_j\times T^k$, and $\dim\mathcal{K}$ measures the dimension of the part of the mixed class aligned with $T^k$.
In such settings, the invariant $r^{\sharp}$ acquires a structural meaning that goes beyond the bound itself: it detects how deeply the torsion coupling penetrates the de Rham irreducible geometry of the factors.
\end{remark}

\begin{remark}[Geometric Significance of Invariance]
\label{rmk:invariance}

The invariance of the lower bound $r^\sharp$ is a non‑trivial property that holds under the Assumption \ref{def:admissible}.  In those cases, the bound depends only on:
\begin{itemize}
    \item[(i)] the mixed class $[\omega]_{\mathrm{mixed}}$,
    \item[(ii)] the parallel‑form strata $\mathcal{P}_k(M_i)$ of the factor metrics.
\end{itemize}
Consequently, $r^\sharp$ is invariant under any deformation of $g_1$ or $g_2$ that 
preserves these strata.  This stability is relevant for manifolds with special holonomy 
where $\mathcal{P}_k(M_i)$ are non‑trivial.
\end{remark}

\begin{remark}[Obstruction to product reducibility and quantified off‑diagonal holonomy]
\label{rmk:forced-irred}
Since \(\nabla^{C}\) is metric, any non‑trivial off‑diagonal holonomy implies that the product splitting \(T_{p}M=V_1\oplus V_2\) is not invariant under \(\mathfrak{hol}_{p}(\nabla^{C})\).  
Theorem \ref{thm:main} quantifies this obstruction to product reducibility in terms of the mixed tensor rank, modulo the obstruction space
\[
\mathcal{K}=\ker\Psi_p\oplus\ker\widetilde{\Psi}_p,
\]
where, in the notation of the theorem,
\[
\ker\Psi_p=\mathcal{V}_{2,1}\cap\bigl(\mathcal P_2(M_1)\otimes \mathcal P_1(M_2)\bigr),\qquad
\ker\widetilde{\Psi}_p=\mathcal{V}_{1,2}\cap\bigl(\mathcal P_1(M_1)\otimes \mathcal P_2(M_2)\bigr).
\]

Specifically, the dimension of the off‑diagonal holonomy subspace satisfies
\[
\dim\mathfrak{hol}_{p}^{\operatorname{off}}(\nabla^{C})\geq 
r^\sharp\;:=\;\operatorname{rank}_{\mathbb{R}}([\omega]_{\mathrm{mixed}})-\dim\mathcal{K},
\]
and the bound certifies that $V_1\oplus V_2$ is not holonomy-invariant when $r^\sharp > 0$, i.e. when $\operatorname{rank}_\mathbb{R}([\omega]_{\mathrm{mixed}}) > \dim\mathcal{K}$; when $r^\sharp = 0$ the inequality is trivially satisfied; the geometric content of this regime is discussed in Remark \ref{rem:deRham}.

In other words, a non‑zero value of the corrected rank
\(r^{\sharp}\) guarantees that the holonomy representation cannot split according to the product decomposition \(T_pM=V_1\oplus V_2\); the connection is geometrically ``entangled'' across the two factors.
\end{remark}

\begin{remark}[Analytic approximation as a tool of proof]
\label{non-analytic}
The proof of the density statement in Lemma \ref{lem:maximal-rank-open} makes essential use of real-analytic approximations of the smooth metric. Real-analytic metrics are dense in the $C^\infty$ topology (Greene–Jacobowitz \cite{GreeneJacobowitz})), and the relevant quantities – the dimension of the space of parallel forms, the harmonic representatives of fixed cohomology classes, and the covariant derivatives of those harmonic forms – depend continuously on the metric in the appropriate topologies. The non-emptiness conclusion of Lemma \ref{lem:maximal-rank-open}, which is all that is required in the proof of Theorem \ref{thm:main}, holds for arbitrary smooth metrics without recourse to analytic approximation.

In Lemma \ref{lem:maximal-rank-open} the reduction to an analytic metric allows us to apply the identity principle for analytic functions: if an analytic section of a vector bundle vanishes on an open set, it vanishes identically.  This principle is used to show that a relation among the covariant derivatives of harmonic forms on an open set forces the existence of a parallel combination of the forms, contradicting the hypothesis of linear independence modulo parallel forms. For a general smooth metric, the non-emptiness of $\mathcal{U}$ is immediate from the non-parallelism hypothesis; the density conclusion is established only under the additional assumption that $g_1$ is real-analytic.

In Theorem \ref{thm:main}, only the non-emptiness of $\mathcal{U}$ is needed: in the admissible setting the kernel identification (Claim 5.2.1) reduces to the non-vanishing of $\nabla^{g_2}\beta_1$ on a non-empty open subset of $M_2$
(Lemma \ref{lem:maximal-rank-open}), and the purely off-diagonal holonomy element is constructed on $\mathcal{W}^{2,1}$ independently of that subset.   

Thus analyticity serves as a tool to establish the stronger density statement in Lemma \ref{lem:maximal-rank-open}; the final statements of Theorem \ref{thm:main} are valid for arbitrary smooth product metrics and do not require any analyticity assumption on the metric or the torsion.

\end{remark}

\begin{remark}[Dichotomy of the obstruction]
\label{rmk:dichotomy}
For admissible product manifolds, the obstruction space $\mathcal{K}$ is controlled by a single topological datum: the vanishing or non-vanishing of $\mathcal{P}_1(M_i)$.  Since $\mathcal{P}_2(M_i)=\mathbb{R}\cdot\operatorname{vol}_{M_i}$ for every compact oriented surface, the obstruction kernels reduce to
\[
\ker\Psi_p\neq\{0\}\;\Longleftrightarrow\;\beta\in\mathcal{P}_1(M_2),
\qquad
\ker\widetilde{\Psi}_p\neq\{0\}\;\Longleftrightarrow\;\alpha\in\mathcal{P}_1(M_1),
\]
where $\alpha$, $\beta$ are the harmonic forms appearing in the minimal decomposition of $[\omega]_{\mathrm{mixed}}$. A compact oriented surface $F$ satisfies $\mathcal{P}_1(F)=\{0\}$ if and only if one of the following holds:
$g(F)=0$ (since $H^1(S^2)=\{0\}$); $g(F)\geq 2$ (by the Poincar\'e--Hopf theorem: $\chi(F)\neq 0$ forbids nowhere-vanishing vector fields, hence parallel $1$-forms);
or $g(F)=1$ and the metric is not flat (a non-zero parallel $1$-form on a surface forces the Gaussian curvature to vanish identically).
In the context of Theorem \ref{thm:main}, the case $g(F)=0$ is vacuous: a non-trivial mixed class in $H^2(M_1)\otimes H^1(M_2)$ requires $H^1(M_2)\neq\{0\}$, hence $g(M_2)\geq 1$.

This yields a sharp dichotomy.  When at least one factor has genus $g\geq 2$, the obstruction space $\mathcal{K}$ is at most one-dimensional regardless of the metric, and $r^\sharp=1$ whenever the mixed class $[\omega]_{\mathrm{mixed}}$ is non-trivial in the corresponding K\"unneth component; in particular, on $\Sigma_{g_1}\times\Sigma_{g_2}$ with $g_1,g_2\geq 2$, the bound $r^\sharp=1$ is \emph{purely topological} and forces the product splitting $V_1\oplus V_2$ to fail to be holonomy-invariant for every cohomologically calibrated connection with torsion representing a non-trivial mixed class, independently of the choice of product metric.  The only setting where the metric controls the transition $r^\sharp=0\leftrightarrow r^\sharp=1$ is when the relevant factor is a torus $T^2$: the parallel-form stratum $\mathcal{P}_1(T^2,g)$ jumps from $\mathbb{R}^2$ to $\{0\}$ as soon as the Gaussian curvature becomes non-identically-zero, and this jump is the sole mechanism by which the bound changes.

Finally, we note that under the pure bidegree hypothesis, exactly one of the two K\"unneth components is active, so $r-\dim\mathcal{K}\leq 1$. We conjecture (Conjecture \ref{conj:mixed}) that the theorem extends to torsion forms of mixed bidegree; this is discussed further in Section \ref{sec:conclusion}.
\end{remark}

\section{Illustrative Examples}
\label{sec:examples}

\subsection*{The case $\Sigma_g \times T^2$ (Example of a non-trivial lower bound)}

Fix $g \geq 2$ and endow $\Sigma_g$ with any hyperbolic metric and
$T^2$ with the flat metric.  The relevant harmonic spaces are
\[
\mathcal{H}^{1}(\Sigma_g)\cong\mathbb{R}^{2g},\quad
\mathcal{H}^{2}(\Sigma_g)
  =\operatorname{span}\{\operatorname{vol}_{\Sigma}\},
\]
\[
\mathcal{H}^{1}(T^{2})
  =\operatorname{span}\{dx,dy\},\quad
\mathcal{H}^{2}(T^{2})
  =\operatorname{span}\{\operatorname{vol}_{T^{2}}\}.
\]
Choose a non-zero harmonic $1$-form
$\alpha\in\mathcal{H}^{1}(\Sigma_g)$
(possible because $b_{1}(\Sigma_g)=2g>0$).
For the hyperbolic metric on $\Sigma_g$,
\[
\mathcal{P}_{1}(\Sigma_g)=\{0\},\qquad
\mathcal{P}_{2}(\Sigma_g)
  =\operatorname{span}\{\operatorname{vol}_{\Sigma}\},
\]
since the volume form is parallel, while no harmonic $1$-form on a
hyperbolic surface is parallel (the holonomy of a hyperbolic metric
acts irreducibly on each fibre of $T^{*}\Sigma_g$, so the only
parallel sections of $\Lambda^{k}T^{*}\Sigma_g$ are proportional
to the volume form for $k=2$ and to the zero section for $k=1$).
For the flat metric on $T^{2}$ all harmonic forms are parallel:
\[
\mathcal{P}_{1}(T^{2})=\mathcal{H}^{1}(T^{2})\cong\mathbb{R}^{2},
\qquad
\mathcal{P}_{2}(T^{2})=\mathcal{H}^{2}(T^{2})\cong\mathbb{R}.
\]

Consider the cohomology class of pure bidegree $(1,2)$:
\[
[\omega]
=\alpha\otimes\operatorname{vol}_{T^{2}}
\;\in\;
\mathcal{H}^{1}(\Sigma_g)\otimes\mathcal{H}^{2}(T^{2})
\;\subset\;
H^{3}(\Sigma_g\times T^{2}).
\]
Since $[\omega]$ lies in a single K\"unneth summand, the pure bidegree
hypothesis of Theorem \ref{thm:main} is satisfied for any metric
connection whose torsion $T$ represents $[\omega]$ and lies in
$\Gamma(V_1^{*}\otimes\Lambda^{2}V_2^{*})$.
The minimal tensor ranks are
\[
r_{2,1}=0,\qquad r_{1,2}=1,\qquad
\operatorname{rank}_{\mathbb{R}}
  \bigl([\omega]_{\mathrm{mixed}}\bigr)=1.
\]

The obstruction kernel for the unique active component is
\[
\ker\widetilde{\Psi}
=\mathcal{V}_{1,2}\cap
  \bigl(\mathcal{P}_{1}(\Sigma_g)
        \otimes\mathcal{P}_{2}(T^{2})\bigr)
=\{0\},
\]
because $\mathcal{P}_{1}(\Sigma_g)=\{0\}$:
the $1$-form $\alpha$ is harmonic but not parallel, so
the simple tensor $\alpha\otimes\operatorname{vol}_{T^{2}}$
does not belong to the parallel stratum.
Hence $\dim\mathcal{K}=0$ and the corrected rank is
\[
r^{\sharp}
=\operatorname{rank}_{\mathbb{R}}
  \bigl([\omega]_{\mathrm{mixed}}\bigr)
  -\dim\mathcal{K}
=1-0=1.
\]

Theorem \ref{thm:main} therefore yields the non-trivial lower bound
\[
\dim\mathfrak{hol}^{\mathrm{off}}_{p}(\nabla^{C})\;\geq\;1.
\]
In particular, for any cohomologically calibrated metric connection with torsion of pure bidegree $(1,2)$ representing $[\omega]$, the product splitting of $\Sigma_g\times T^{2}$ is not invariant under the holonomy (the connection does not reduce to a product connection).

\subsection*{The case $T^2 \times T^2$ 
(Example of a trivial lower bound)}

Let $M = T^2 \times T^2$ endowed with the flat product 
metric $g = g_1 \oplus g_2$.  Fix harmonic bases
\[
  H^1(T^2) = \operatorname{span}\{dx, dy\}, \quad 
  H^2(T^2) = \operatorname{span}\{\operatorname{vol}\}.
\]
For the flat metric we have
\[
  \mathcal{P}_1(T^2) = H^1(T^2) \cong \mathbb{R}^2, 
  \qquad
  \mathcal{P}_2(T^2) = H^2(T^2) \cong \mathbb{R}.
\]
The mixed K\"unneth components are
\[
  H^2(T^2) \otimes H^1(T^2) \cong \mathbb{R}^2, \quad 
  H^1(T^2) \otimes H^2(T^2) \cong \mathbb{R}^2.
\]
Consider the cohomology class
\[
  [\omega] = \underbrace{\operatorname{vol} \otimes dx}
  _{(2,1)} 
  + \underbrace{dy \otimes \operatorname{vol}}_{(1,2)}.
\]
Its minimal tensor-rank decomposition is
\[
  \pi_{2,1}([\omega]) = \operatorname{vol} \otimes dx 
  \quad (\text{rank } 1), \qquad
  \pi_{1,2}([\omega]) = dy \otimes \operatorname{vol} 
  \quad (\text{rank } 1).
\]
Hence
\[
  \operatorname{rank}_{\mathbb{R}}
  ([\omega]_{\text{mixed}}) = 1 + 1 = 2.
\]
Since $\operatorname{vol} \in \mathcal{P}_2$ and 
$dx, dy \in \mathcal{P}_1$, both simple tensors lie in 
the obstruction spaces
\[
  \ker\Psi_p = \operatorname{span}
  \{\operatorname{vol} \otimes dx\}, \qquad
  \ker\tilde{\Psi}_p = \operatorname{span}
  \{dy \otimes \operatorname{vol}\},
\]
so $\dim\mathcal{K} = 2$ and
\[
  r^{\sharp} = \operatorname{rank}_{\mathbb{R}}
  \bigl([\omega]_{\mathrm{mixed}}\bigr) 
  - \dim\mathcal{K} = 2-2 = 0.
\]
We note that the class $[\omega]$ has non-pure bidegree, 
since both its $(2,1)$ and $(1,2)$ K\"unneth components 
are non-trivial, and therefore falls outside the strict 
hypotheses of Theorem \ref{thm:main}.  However, when 
$r^{\sharp}=0$ the conclusion 
$\dim\mathfrak{hol}_{p}^{\mathrm{off}}(\nabla^{C})
\geq 0$ is universally true and requires no hypothesis 
on the bidegree of $T$.

This example illustrates the dichotomy of 
Remark \ref{rmk:dichotomy}: when both factors are flat 
tori, the parallel-form strata are maximal 
($\mathcal{P}_k(T^2) = H^k(T^2)$ for all $k$) and the 
obstruction kernel absorbs the entire mixed class, 
yielding a trivially saturated bound; the geometric interpretation of this regime is discussed in Remark \ref{rem:deRham}.

\section{Conclusion}
\label{sec:conclusion}

We have shown that for \emph{admissible} product manifolds \(M = M_1 \times M_2\), specifically, products of compact oriented surfaces, and for any cohomologically calibrated metric connection \(\nabla^C\) with torsion of pure bidegree, the off-diagonal holonomy dimension admits a quantified lower bound on a non-empty open subset $\mathcal{V} \subset M$
\[
\dim\bigl(\mathfrak{hol}_p^{\mathrm{off}}(\nabla^{C})\bigr)\;\geq\;
r^\sharp\;:=\;\operatorname{rank}_{\mathbb{R}}\bigl([\omega]_{\mathrm{mixed}}\bigr)-\dim\mathcal{K},
\]
where \(r=\operatorname{rank}_{\mathbb{R}}([\omega]_{\operatorname{mixed}})\) and \(\mathcal{K}=\ker\Psi_p\oplus\ker\widetilde{\Psi}_p\) is the sum of the kernels of the linear obstruction maps. The bound is metric‑invariant under deformations preserving the parallel‑form strata \(\mathcal{P}_k(M_i)\). This provides the first \emph{computable} obstruction to the reduction of the holonomy along the product splitting on product manifolds within the admissible class.
As shown in Remark \ref{rmk:dichotomy}, the behaviour of $r^\sharp$ across the admissible class is completely determined: $r^\sharp=1$ is forced by topology alone whenever at least one factor has genus $g\geq 2$ and the mixed class is non-trivial, while the transition $r^\sharp=0\leftrightarrow 1$ occurs only when a torus factor changes from flat to non-flat.

The examples in Section \ref{sec:examples} illustrate both non‑trivial  and trivial cases, while Remark \ref{rem:counterexample} shows by explicit construction that the bound fails to be well‑defined without the admissibility hypothesis, establishing sharpness of our assumptions.

When the parallel-form strata absorb the entire mixed class, $r^\sharp = 0$: no de-Rham--irreducible factor prevents the absorption, and transverse non-reducibility (Remark \ref{rem:deRham}) is absent. However, $r^\sharp = 0$ does not imply that the off-diagonal holonomy vanishes: as observed in Remark \ref{rem:harmonic}, the quadratic torsion term $\frac{1}{4}Q_T$ can generate non-trivial off-diagonal curvature even when the harmonic projection $\omega_h$ is entirely parallel.
This reflects a deliberate feature of the construction:
$r^\sharp$ is defined through the harmonic projection $\omega_h$ alone, because for a fixed cohomology class $[\omega]$ infinitely many torsion forms represent $[\omega]$, so any quantity depending on the full torsion would fail to be an invariant of the class.
This invariance ensures that $r^\sharp$ is computable from the
topology of $M$ and the parallel-form strata of the factor
metrics, without knowledge of the specific torsion form.

Strict inequality $\dim\mathfrak{hol}^{\mathrm{off}}_p(\nabla^C) > r^\sharp$ remains an open question.

\begin{corollary}[Extended admissible setting]
\label{cor:extended}
Theorem \ref{thm:main} extends to products 
$M_1 \times M_2$ where only the factor contributing 
the $2$-form to the torsion bidegree is required to 
be a compact oriented surface, while the other factor 
may be a compact oriented Riemannian manifold of 
arbitrary dimension.  Specifically:
\begin{enumerate}
  \item[\textup{(i)}] For torsion of pure bidegree 
    $(2,1)$, $M_1$ is required to be a compact oriented 
    surface; $M_2$ may have arbitrary dimension.
  \item[\textup{(ii)}] For torsion of pure bidegree 
    $(1,2)$, $M_2$ is required to be a compact oriented 
    surface; $M_1$ may have arbitrary dimension.
\end{enumerate}
When both factors are surfaces, both bidegree types are covered simultaneously, recovering the full  admissible setting of 
Definition \ref{def:admissible}.
\end{corollary}

\begin{proof}
The proof of Theorem \ref{thm:main} for pure bidegree 
$(2,1)$ uses the hypothesis $\dim M_1 = 2$ but not 
$\dim M_2 = 2$; symmetrically, for pure bidegree 
$(1,2)$ the proof uses $\dim M_2 = 2$ but not 
$\dim M_1 = 2$.  No modification of \emph{Claim 5.2.1} is required: its kernel computation uses only that $\alpha_1\in\mathcal{H}^2(M_1)=\mathbb{R}\cdot\operatorname{vol}_{M_1}$ is parallel --- which holds because $M_1$ is a surface --- together with Lemma \ref{lem:maximal-rank-open} applied to $\beta_1$ on $M_2$, valid for $M_2$ of
arbitrary dimension. Symmetrically, for pure bidegree $(1,2)$ one uses that $\widetilde\beta_1\in\mathcal{H}^2(M_2)=\mathbb{R}\cdot\operatorname{vol}_{M_2}$ is parallel, $M_2$ a surface, with $M_1$ of arbitrary dimension. The remainder of the proof is unchanged.
In particular Lemma \ref{lem:invariance} and Definition \ref{def:obstruction} extend verbatim to the $(2,1)$ case in this setting, since they use only $\mathcal{H}^2(M_1)=\mathbb{R}\cdot\operatorname{vol}_{M_1}$; under pure bidegree $(2,1)$ one has $\mathcal{V}_{1,2}=\{0\}$, so $\mathcal{K}=\ker\Psi_p$.
\end{proof}

\subsection*{Dimensional reduction along a circle factor}

In this extended setting, the classification of $\mathcal{P}_1(M_j)$ in Remark \ref{rmk:dichotomy} generalises from genus and Gaussian curvature to the de Rham decomposition theory of the higher-dimensional factor, and the geometric 
interpretation of $r^\sharp$ (Remark \ref{rem:deRham}) 
acquires additional depth: $r^\sharp$ detects whether the 
cohomological coupling penetrates beyond the flat sub-factors 
identified by parallel $1$-forms on $M_j$. A concrete instance arises when $M_2 = N \times S^1$ (respectively $M_1 = N \times S^1$ for bidegree $(1,2)$): the invariant $r^\sharp$ then predicts whether the failure of the connection to reduce along the product splitting survives dimensional reduction along the circle factor. 
When $r^\sharp = 1$, the torsion $3$-form survives restriction to $M_1 \times N$, and the restricted connection still fails to preserve the product splitting of the reduced product; when $r^{\sharp}=0$ with $\beta\in\mathbb{R}\,d\theta$, the restricted class is trivial, $[\iota_\theta^{*}T]=0$: the harmonic projection of the restricted torsion is no longer a non-trivial mixed class, so the reduced connection is not
cohomologically calibrated and the mechanism detected by $r^{\sharp}$ is absent on $M_1\times N$.  Since the restricted class vanishes it admits the zero representative, whose associated torsion-free connection has holonomy preserving
the product splitting of the reduced product; for a general representative the restricted torsion need not vanish, and by \emph{Case 1} it may still generate off-diagonal holonomy.

The first of these assertions is made precise in Proposition \ref{prop:reduction} below.

\begin{proposition}[Reduction along a circle factor]\label{prop:reduction}
Let $M_1$ be a compact connected oriented surface, $M_2=N\times S^1$ with a product metric, $N$ compact connected oriented, and let $\iota_\theta:M_1\times N\hookrightarrow M$ be the slice at $\theta$.  Let $\nabla^C$ be cohomologically calibrated with torsion of pure bidegree $(2,1)$,
$\omega_h=\operatorname{vol}_{M_1}\wedge\beta$.  If $r^{\sharp}=1$, then $\iota_\theta^{*}T\neq0$ for some $\theta$, and the metric connection on $M_1\times N$ with skew torsion $\iota_\theta^{*}T$ satisfies $\dim\mathfrak{hol}^{\mathrm{off}}\geq1$ on a non-empty open set; in particular its holonomy does not preserve the product splitting.
\end{proposition}

\begin{proof}
Write $T=\operatorname{vol}_{M_1}\wedge\sigma$ and $\sigma=\sigma_N+u\,d\theta$.
As $\iota_\theta^{*}d\theta=0$, $\iota_\theta^{*}T=0$ for all $\theta$ iff $\sigma_N\equiv0$, i.e. $\sigma=u\,d\theta$.  In that case the pointwise orthogonality $T^{*}N\perp\mathbb{R}\,d\theta$ gives $\langle T,\operatorname{vol}_{M_1}\wedge\beta'\rangle_{L^2}=0$ for every $\beta'\in\mathcal{H}^{1}(N)$; since $\operatorname{vol}_{M_1}\wedge\beta'$ is harmonic, this equals $\langle\omega_h,\operatorname{vol}_{M_1}\wedge\beta'\rangle_{L^2}$, and taking $\beta'$ equal to the $\mathcal{H}^{1}(N)$-component of $\beta$ gives that this component vanishes; since
$\mathcal{H}^{1}(M_2)=\mathcal{H}^{1}(N)\oplus\mathbb{R}\,d\theta$ and the
harmonic projection preserves the bidegree, $\beta\in\mathbb{R}\,d\theta
\subseteq\mathcal{P}_1(M_2)$, whence $r^{\sharp}=0$.  The first assertion is the contrapositive.  Since $\iota_\theta^{*}T\in\Gamma(\Lambda^2V_1^{*}\otimes T^{*}N)$
has pure bidegree $(2,1)$ and $M_1$ is a surface, \emph{Case 1} of the proof of Theorem \ref{thm:main} applies in the extended setting of Corollary \ref{cor:extended} at every point where $\iota_\theta^{*}T\neq0$.
\end{proof}

\begin{conjecture}\label{conj:mixed}
Theorem \ref{thm:main} extends to torsion forms of mixed bidegree, i.e.\ the pure bidegree hypothesis can be removed.
\end{conjecture}

Under mixed bidegree both K\"unneth components may be simultaneously active, so $r-\dim\mathcal{K}$ can reach the value $2$ in the admissible setting. Establishing the conjecture would require new techniques for the overlap region $\mathcal{W}^{2,1}\cap\mathcal{W}^{1,2}$, where both components of the torsion are simultaneously non-vanishing.  A systematic investigation of this extension remains an interesting open problem.

A further direction is to explore extensions beyond the admissible setting. Imposing \emph{Kruskal‑rank constraints} on harmonic representatives could yield intrinsic mixed factor spaces even when \(b_2(M_i) > 1\), potentially treating Calabi–Yau threefolds or \(G_2 \times S^1\) backgrounds. Such extensions would require new algebraic tools and are left for future investigation.

Another further interesting direction for future work could be to explore the implications of these results in the context of deformations of complex structures, a field pioneered by the works of Kodaira \cite{K (1958)} and Kodaira-Spencer \cite{K-S (1960)}.

\medskip
\noindent\textbf{Compactness and Hodge theory.}
Compactness ensures the existence and uniqueness of harmonic representatives \cite{W (1983)}, and the vanishing of boundary terms in integrations by parts. Extending the results to non-compact manifolds would require an $L^2$ Hodge framework or alternative analytic hypotheses.

\end{document}